\documentclass{article}
\usepackage{graphicx} 

\title{The maximum size of $3$-distance binary codes}

\usepackage{amsmath, amssymb, amsthm, url, tikz, verbatim}

\newtheorem{theorem}{Theorem}[section]

\newtheorem{lemma}[theorem]{Lemma}

\newtheorem{question}[theorem]{Question}

\newtheorem{conjecture}[theorem]{Conjecture}

\DeclareMathOperator{\Tr}{Tr}

\DeclareMathOperator{\im}{im}
\DeclareMathOperator{\spa}{span}

\makeatletter
\NewCommandCopy\@@pmod\pmod
\DeclareRobustCommand{\pmod}{\@ifstar\@pmods\@@pmod}
\def\@pmods#1{\mkern4mu({\operator@font mod}\mkern 6mu#1)}
\makeatother

\usetikzlibrary{calc}
\usetikzlibrary{patterns}
\usetikzlibrary{decorations.markings}
\usetikzlibrary{arrows,shapes.geometric,%
	decorations.pathreplacing,shapes,shadows}

\date{}

\author{Valjakas Djaljapayan\thanks{
Department of Mathematics, National Central University,
Taoyuan City, 320, Taiwan. Valjakas Djaljapayan and Wei-Hsuan Yu are supported by MOST
under Grant No.~115-2115-M-008-007-MY2. Email: {\tt 115281004@cc.ncu.edu.tw; whyu@math.ncu.edu.tw}.
}
~~~~
Hao Huang\thanks{Department of Mathematics, National University of Singapore, 119076, Singapore. Research supported by a start-up grant at NUS and an MOE Academic Research Fund (AcRF) Tier 1 grant A-8003627. Email: {\tt huanghao@nus.edu.sg}.
}
~~~~
Wei-Hsuan Yu\footnotemark[1]
}

\begin{document}

\maketitle

\begin{abstract}
In this paper, we show that for every positive integer $n \ge 24$, the maximum size of a family $\mathcal{F}$ of vectors from $\{-1,1\}^n$ with at most three distinct Hamming distances is equal to $\binom{n}{3}+n$. The threshold $n=24$ is best possible.
\end{abstract}
\section{Introduction}

Let $(Y,\operatorname{dist})$ be a metric space. A finite subset
$X\subseteq Y$ is called an \emph{$s$-distance set} if the set
\[
D(X):=\{\operatorname{dist}(x,y):x,y\in X,\ x\ne y\}
\]
has cardinality at most $s$. A basic problem is to determine, for a
given space $Y$ and a given integer $s$, the largest possible
cardinality of an $s$-distance set in $Y$. As the underlying space varies, this problem connects discrete geometry, algebraic
combinatorics, spectral graph theory, and coding theory. 

The geometric origins of this problem precede the modern terminology.
Early work includes Kelly's 1947 solution of Erd\H{o}s's problem on isosceles
configurations \cite{EK47} and Croft's study of finite configurations in
three-dimensional space \cite{Croft62}. These are low-dimensional antecedents;
the systematic study of Euclidean two-distance sets in arbitrary dimension
was developed by Einhorn and Schoenberg in their two papers of 1966
\cite{ES66a,ES66b}. Their work connected distance configurations with linear
algebra and graph representations and established basic finiteness and
classification results. Already the simplest examples illustrate the
extremal question: a regular simplex gives $n+1$ equidistant points in
$\mathbb R^n$, while a regular pentagon is a two-distance set of five points
in the plane.

For an $s$-distance subset $X$ of $\mathbb{R}^{n}$ with the usual Euclidean distance, Bannai, Bannai, and Stanton
\cite{BannaiBannaiStanton} proved the general bound
\[
|X|\le \binom{n+s}{s}.
\]
The study of equality cases and improvements to sharper bounds has led to
connections with spherical designs, association schemes, and the
spectra of graphs; see, for example,
\cite{JiangTidorYaoZhangZhao,LarmanRogersSeidel,Nozaki,NozakiSuda}.
Further work on Euclidean and spherical few-distance sets includes
\cite{Lisonek,MusinNozaki,OstergardSzollosi,YuLiu}, while
\cite{GlazyrinYu,HegedusRonyai,PetrovPohoata} illustrate the role of
polynomial and spectral methods in more general settings.

In this paper we study the analogous extremal problem in the binary Hamming space $\{-1,1\}^n$. For two vectors $x, y \in \{-1,1\}^n$, we denote by $d(x,y)$  the Hamming distance between $x$ and $y$, i.e. the number of coordinates $i \in [n]:=\{1,\ldots,n\}$ such that $x_i \neq y_i$. We say that $\mathcal{F}$ is an $s$-distance set if $\{d(x,y): x, y \in \mathcal{F}, x\neq y\}$ has size at most $s$. One can easily use the polynomial method to show that for every $s$-distance set $\mathcal{F} \subset \{-1, 1\}^n$, 
\begin{align}\label{ineq_trivial}
    |\mathcal{F}| \le \sum_{i=0}^s \binom{n}{i}.
\end{align}
This upper bound can be attained by some well-known constructions. For example, starting from any $(n+1)\times(n+1)$ Hadamard
matrix with entries from $\{-1,1\}$ and an all-one first column,
we remove its first column and take the remaining row vectors.
This gives a $1$-distance set
$\mathcal{F}\subset\{-1,1\}^n$ with $|\mathcal{F}|=n+1$. For $n=2s+1$, taking all the vectors with an even number of $1$'s gives a family with $|\mathcal{F}|=2^{2s}=\sum_{i=0}^s \binom{n}{i}$, and such $\mathcal{F}$ is an $s$-distance set since the possible Hamming distances are $\{2, 4, \cdots, 2s\}$. For $s=3$, the dual binary Golay code gives a $3$-distance set in $\{-1,1\}^{23}$ of size $2048=1+23+\binom{23}{2}+\binom{23}{3}$.

However, in all these constructions except for the Hadamard ones,  the dimension $n$ turns out to be either a fixed constant (as in the dual Golay construction), or a function of the number of distances $s$. For $s \ge 2$, there are no known $s$-distance families of vectors $\mathcal{F} \subset \{-1,1\}^n$ attaining this upper bound as $n \rightarrow \infty$. This naturally leads to the following folklore conjecture, where only roughly half of the binomial coefficients remain in the sum in the upper bound.

\begin{conjecture}
\textnormal{(\cite{BargMusin,DongGaoLiuOuyangZhou,MusinNozaki})}
\label{conj_main}
For every integer $s\ge2$, there exists a positive integer $M_s$
such that, whenever $n\ge M_s$, every $s$-distance set
$\mathcal F\subseteq\{-1,1\}^n$ satisfies
\[
|\mathcal F|
\le
\sum_{\substack{0\le i\le s\\ i\equiv s\pmod* 2}}\binom{n}{i}.
\]
\end{conjecture}

For $n\ge2s$, the conjectured bound is attained by the family
of vectors whose numbers of $1$-coordinates are at most $s$
and congruent to $s$ modulo $2$.
This family has distance set $\{2,4,\ldots,2s\}$; see also
\cite[Eq.~(3.1)]{MusinNozaki}.
The restriction $s\ge2$ is necessary, since Hadamard matrices
give equidistant sets of $n+1$ vectors in $\{-1,1\}^n$
for arbitrarily large $n$.
Moreover, when $n=2s+1$, the full even parity class has
distance set $\{2,4,\ldots,2s\}$ and size $2^{2s}$,
exceeding the conjectured bound.
Hence any threshold $M_s$ must satisfy $M_s\ge2s+2$.

Barg and Musin \cite{BargMusin} obtained general upper bounds
and exact results for binary codes with few distances.
Their Theorem~11 shows that, if the prescribed distances
satisfy $d_1+\cdots+d_s\le sn/2$, then
\[
|\mathcal F|
\le
\binom ns+\sum_{i=0}^{s-2}\binom ni.
\]
Combining harmonic bounds with linear programming, they also
determined the maximum size of binary $3$-distance sets to be
$n+\binom n3$ for $8\le n\le22$ and $n=24$, and the maximum
size of binary $4$-distance sets to be
$1+\binom n2+\binom n4$ for $10\le n\le33$
\cite[Proposition~13]{BargMusin}.

Musin and Nozaki \cite{MusinNozaki} extended these exact results
by combining harmonic and linear programming bounds with
integrality conditions from a generalization of the
Larman--Rogers--Seidel theorem.
Their Theorem~3.2 establishes the same three-distance maximum
$n+\binom n3$ for the additional lengths
$24\le n\le33$ and $n=36,37,44$.
For four distances, they extended the range of optimality
of $1+\binom n2+\binom n4$ to all lengths $10\le n\le47$.
Thus these two papers established the conjectured maximum
for substantial finite ranges in both the three-distance
and four-distance cases.

Barg, Glazyrin, Kao, Lai, Tseng, and Yu
\cite{BargGlazyrinKaoLaiTsengYu}
subsequently resolved the case $s=2$ of
Conjecture~\ref{conj_main}, proving that the maximum size
is $1+\binom n2$ for every $n\ge6$.
Thus the optimal threshold is $M_2=6$.
For three distances, they added the cases $n=34,35$
\cite[Proposition~6.1]{BargGlazyrinKaoLaiTsengYu},
completing the known interval $24\le n\le37$.

For a finite set $D$ of positive integers, let $A(n,D)$ denote
the maximum size of a family $\mathcal F\subseteq\{-1,1\}^n$
whose pairwise distances belong to $D$.
Dong, Gao, Liu, Ouyang, and Zhou
\cite[Theorem~1.4]{DongGaoLiuOuyangZhou}
proved that, for every fixed $s$-element set
$D\ne\{2,4,\ldots,2s\}$,
\[
A(n,D)
\le
\left(\frac{s}{s+1}+o_D(1)\right)\binom ns
\qquad\text{as }n\to\infty.
\]
They also established
\[
A(n,\{2,4,\ldots,2s\})
=
\sum_{\substack{0\le i\le s\\ i\equiv s\pmod* 2}}
\binom{n}{i}
\qquad\text{for every }n\ge2s+2.
\]
Thus the conjectured inequality holds for each fixed distance
set once $n$ is sufficiently large in terms of that set.
Conjecture~\ref{conj_main}, however, requires a threshold
depending only on $s$, allowing the distances to vary with $n$.

In this paper, we resolve the case $s=3$ of
Conjecture~\ref{conj_main}, uniformly over all choices of
distances, for every $n\ge24$.
The threshold is optimal, since the dual binary Golay code
of length $23$ has $2048>23+\binom{23}{3}$ codewords
\cite[Remark~3.3]{MusinNozaki}.
Our main result is the following.

\begin{theorem}\label{thm_main}
Suppose $0<d_1<d_2<d_3 \le n$, and $\mathcal{F}$ is a family of vectors from $\{-1,1\}^n$ such that the Hamming distance between any of its two distinct vectors is from $\{d_1, d_2, d_3\}$. Then when $n \ge 24$, we always have
$$|\mathcal{F}| \le \binom{n}{3}+n.$$
\end{theorem}

The bound is attained by the family of all vectors with exactly one
or exactly three $1$-coordinates. Consequently, for every $n\ge 24$,
the maximum size of a binary $3$-distance set is exactly
\[
n+\binom{n}{3}.
\]
The threshold is best possible: the dual binary Golay code gives a
$3$-distance set of size $2048$ when $n=23$, whereas
\[
n+\binom{n}{3}=23+\binom{23}{3}=1794<2048.
\]
Therefore the optimal threshold in the $s=3$ case of
Conjecture~\ref{conj_main} is $M_3=24$.

The rest of the paper is organized as follows. In the next section, we state and give short proofs of two useful bounds in extremal set theory and coding theory that we will repeatedly use throughout the paper: the inertia (harmonic) bound and the ratio (linear programming) bound. We also include a simple fact that in order to prove Theorem \ref{thm_main}, it suffices to only consider the case when all the distances are even. The following four sections cover several different ranges: Section \ref{sec_outside} studies the relatively easy case when $d_1+d_2+d_3$ is either at most $\frac{3n}{2}$ or strictly greater than $\frac{3n+6}{2}$. The proof of the first case uses the polynomial method, and the second case uses the inertia bound. It remains to investigate the exceptional interval $d_1+d_2+d_3 \in [\frac{3n+1}{2}, \frac{3n+6}{2}]$. Section \ref{sec_even} resolves all even $n$ cases by a combination of the inertia bound and the ratio bound, as well as some trace analysis. Section \ref{sec_odd} solves all the odd $n$ cases except for a difficult case when $d_1+d_2+d_3=\frac{3n+3}{2}$ and a certain function $c_1$ of $n$ and the distances $d_i$ is strictly less than $6$. This difficult case is separately addressed and proved in Section \ref{sec_difficult} based on a series of linear algebra arguments. Finally, Section~\ref{sec_concluding} presents concluding remarks
and open problems concerning the classification of extremal families,
uniform gaps, and extensions to general $s$-distance sets.\\

\noindent{\bf Disclosure of AI use: } The authors used ChatGPT to assist in developing and implementing their ideas. Throughout this process, the authors guided and directed
the use of ChatGPT, reviewed its contributions, and take full
responsibility for the content of the paper, including its mathematical results and references.

\section{The inertia and ratio bounds} \label{sec_bounds}
In this section, we establish two bounds that will be repeatedly used in the proofs. Readers familiar with spectral graph theory will recognize the first as Cvetkovi\'c's inertia bound \cite{Cv71}, which bounds the independence number of a graph by the number of nonnegative eigenvalues of any of its pseudo-adjacency matrices. The second is Hoffman's ratio bound (see \cite{HaemersRatioBound}), which bounds the independence number by the number of vertices multiplied by $\frac{-\lambda_{\min}}{\lambda_{\max}-\lambda_{\min}}$, for any pseudo-adjacency matrix with equal row sum that gives the largest eigenvalue. In the combinatorial geometry and coding theory community, these are commonly known as the harmonic and linear programming bounds, respectively. We include proofs of both results to make the paper self-contained and adopt notation suited to our subsequent applications.

For a subset $T \subset [n]$ and $x \in \mathbb{R}^n$, let $\chi_T(x)=\prod_{i \in T} x_i$.
\begin{lemma}[The inertia/harmonic bound] \label{lem_inertia}
If the matrix $K$ given by $K(x, y)= \sum_T k_T \chi_T(x)\chi_T(y)$, when restricted to $x, y \in \mathcal{F} \subset \{-1, 1\}^n$,  is equal to the matrix $\lambda I$ for some $\lambda>0$, then 
$$|\mathcal{F}| \le |\{T: k_T>0\}|.$$
\end{lemma}
\begin{proof}
We define $v_T=(\chi_T(x))_{x \in \mathcal{F}}$. If there exists a vector $z \in \mathbb{R}^{\mathcal{F}}$  orthogonal to all $v_T$ with $k_T>0$, then 
$$ \lambda \sum_{x \in \mathcal{F}} z_x^2 = \sum_{x, y \in \mathcal{F}} K(x,y) z_x z_y = \sum_{T: k_T \le 0} k_T \langle v_T, z\rangle ^2 \le 0,$$
implying $z=0$. Therefore $\{v_T\}_{T: k_T>0}$ span $\mathbb{R}^{\mathcal{F}}$; a dimension argument immediately gives
$$|\mathcal{F}|=\dim \mathbb{R}^{\mathcal{F}} \le |\{T: k_T>0\}|.$$
\end{proof}

The following ratio bound is also true.
\begin{lemma}[The ratio/LP bound \cite{HaemersRatioBound}] \label{lem_ratio}
Suppose the matrix $K$ given by $K(x, y)= \sum_T k_T \chi_T(x)\chi_T(y)$, when restricted to $x, y \in \mathcal{F} \subset \{-1, 1\}^n$, is equal to the matrix $\lambda I$ for some $\lambda>0$. If, in addition, every $k_T \ge 0$ and $k_{\emptyset}>0$, then 
$$|\mathcal{F}| \le \frac{\lambda}{k_{\emptyset}}.$$
\end{lemma}
\begin{proof}
If $\mathcal F=\varnothing$, the assertion is immediate.
Hence we may assume that $|\mathcal F|>0$. Observe that  
\begin{align*}
\lambda|\mathcal{F}|=\sum_{x \in \mathcal{F}}\sum_{y \in \mathcal{F}} K(x,y)=\sum_T k_T \left(\sum_{x \in \mathcal{F}} \chi_T(x)\right)^2 \ge k_{\emptyset}|\mathcal{F}|^2.
\end{align*}
The last inequality comes from ignoring all the summands with $T \neq \emptyset$ and using $\chi_{\emptyset}=1$. As a consequence, $|\mathcal{F}| \le \frac{\lambda}{k_{\emptyset}}$.
\end{proof}

Our next lemma is a simple observation that  having even just one odd distance would lead to an upper bound quadratic in $n$.
\begin{lemma}
For $n \ge 24$, if $\mathcal{F} \subset \{-1,1\}^n$ is a $3$-distance set with an odd distance, then $|\mathcal{F}| \le n^2+n+2 < \binom{n}{3}+n$.
\end{lemma}
\begin{proof}
We split $\mathcal{F}$ into two subfamilies $\mathcal{F}_0$
and $\mathcal{F}_1$ according to the parity of the number of $1$-coordinates. Then two vectors in the same subfamily have even Hamming distance. Each family $\mathcal{F}_i$ is a $2$-distance set since the odd distance cannot come from two vectors from the same subfamily. Therefore from \eqref{ineq_trivial}, we have $|\mathcal{F}_i| \le \binom{n}{2}+n+1$ for $i=0,1$ and thus 
$$|\mathcal{F}| \le 2 \left(\binom{n}{2}+\binom{n}{1}+1\right)=n^2+n+2<\binom{n}{3}+n$$
for $n \ge 24$.
\end{proof}

\noindent\textbf{Remark.~}
The same argument applies to $s$-distance sets. Suppose that
$s \ge 2$ and $\mathcal{F} \subseteq \{-1,1\}^n$ is an
$s$-distance set with at least one odd distance. Splitting
$\mathcal{F}$ into two subfamilies according to the parity of the
number of $1$-coordinates, we see that each subfamily is an
$(s-1)$-distance set. Therefore, by \eqref{ineq_trivial},
\[
|\mathcal{F}|
\le 2\sum_{i=0}^{s-1}\binom{n}{i}.
\]
In particular, when $n \ge 4s$,
\[
|\mathcal{F}|
< \binom{n}{s}
\le \sum_{\substack{0\le i\le s\\ i\equiv s\pmod* 2}}
\binom{n}{i}.
\]
Indeed, for $n\ge 4s$ and $1\le i\le s-1$, we have $\frac{\binom{n}{i-1}}{\binom{n}{i}}
=\frac{i}{n-i+1}<\frac13.$ Thus
\[
2\sum_{i=0}^{s-1}\binom{n}{i}
<3\binom{n}{s-1}
<\binom{n}{s},
\]
as required.

We may assume that all three prescribed distances occur in
$\mathcal F$. Indeed, if at most two distances occur, then
\eqref{ineq_trivial} gives $|\mathcal F|\le 1+n+\binom{n}{2}<n+\binom{n}{3}$ for $n\geq 24$. If one of the three occurring distances is odd, the preceding lemma
gives the desired bound. Hence, from now on, we may assume that $d_1,d_2,d_3$ are all even.

\section{When the distance sum is outside an exceptional interval}\label{sec_outside}
 For two distinct vectors $x, y \in \mathcal{F} \subset \{-1,1\}^n$, we have $$\langle x, y \rangle=n-2d(x, y) \in \{n-2d_1, n-2d_2, n-2d_3\}.$$
We set $\alpha_i:=n-2d_i$ and define a degree-$3$ single-variable polynomial $P$ as follows:
\begin{align*}
P(z)&=(z-\alpha_1)(z-\alpha_2)(z-\alpha_3)\\
&=z^3-(\alpha_1+\alpha_2+\alpha_3)z^2+(\alpha_1 \alpha_2+\alpha_1 \alpha_3+\alpha_2 \alpha_3)z-\alpha_1 \alpha_2 \alpha_3.
\end{align*}
Whenever $x, y \in \mathcal{F}$ and $x \neq y$, since $\langle x, y\rangle=n-2d(x,y) \in \{\alpha_1, \alpha_2, \alpha_3\}$, we have $P(\langle x, y\rangle)=0$.  On the other hand, for $x=y$, we have $\langle x, y\rangle =n$ and in this case $P(n)=(n-\alpha_1)(n-\alpha_2)(n-\alpha_3)=8d_1d_2d_3$. Let us define the $|\mathcal{F}| \times |\mathcal{F}|$ matrix $M$ with $M_{x, y}=P(\langle x, y\rangle)$ for all $x, y \in \mathcal{F}$. Then clearly $M=8d_1d_2d_3 I$. 

On the other hand, for each $S \subset [n]$, we define a vector $v_S \in \mathbb{R}^{\mathcal{F}}$, such that for every $x=(x_1, \cdots, x_n) \in \mathcal{F}$, $(v_S)_{x}=\chi_S(x)=\prod_{i \in S} x_i$. Using $x_i^2=y_i^2=1$ for every $i$, we have
$$\langle x, y\rangle^2=\left(\sum_{i=1}^n x_i y_i\right)^2=n + 2\sum_{|S|=2} \prod_{i \in S} x_iy_i=n+2 \sum_{|S|=2} (v_S)_{x} (v_S)_{y}.$$
\begin{align*}
\langle x, y\rangle^3&=\left(\sum_{i=1}^n x_i y_i\right)^3=(3n-2) \langle x,y \rangle + 6\sum_{|S|=3} \prod_{i \in S} x_iy_i\\
&=(3n-2) \sum_{|S|=1} (v_S)_{x} (v_S)_{y} + 6 \sum_{|S|=3} (v_S)_{x} (v_S)_{y}
\end{align*}
Thus 
\begin{align*}
M_{x,y}&=P(\langle x, y\rangle)=\langle x, y\rangle^3-(\alpha_1+\alpha_2+\alpha_3)\langle x, y\rangle^2+(\alpha_1\alpha_2+\alpha_1\alpha_3+\alpha_2\alpha_3)\langle x, y\rangle -\alpha_1\alpha_2\alpha_3\\
&=(-\alpha_1\alpha_2\alpha_3-n(\alpha_1+\alpha_2+\alpha_3))+(3n-2+\alpha_1\alpha_2+\alpha_1\alpha_3+\alpha_2\alpha_3)\sum_{|S|=1}(v_S)_x(v_S)_y\\
&~~~~-2(\alpha_1+\alpha_2+\alpha_3) \sum_{|S|=2}(v_S)_x(v_S)_y + 6 \sum_{|S|=3}(v_S)_x(v_S)_y.
\end{align*}
Hence we have the matrix identity
\begin{align}\label{eq_key}
8d_1d_2d_3 I = M= c_0v_{\emptyset} v_{\emptyset}^T + c_1\sum_{|S|=1}v_{S} v_{S}^T + c_2\sum_{|S|=2}v_{S} v_{S}^T + c_3\sum_{|S|=3}v_{S} v_{S}^T
\end{align}
with 
$$c_0=-\alpha_1\alpha_2\alpha_3-n(\alpha_1+\alpha_2+\alpha_3),\qquad c_1=3n-2+(\alpha_1\alpha_2+\alpha_1\alpha_3+\alpha_2\alpha_3),$$
$$c_2=-2(\alpha_1+\alpha_2+\alpha_3),\qquad c_3=6.$$

We will repeatedly use \eqref{eq_key} throughout this paper. Our first theorem settles the case when the sum of $d_i$ is ``small''.
\begin{theorem} \label{thm_n_large}
For $n\ge 24$, if $\mathcal{F} \subset \{-1,1\}^n$ is a $3$-distance set with distances from $\{d_1, d_2, d_3\}$ and $d_1+d_2+d_3 \le \frac{3n}{2}$, then 
$$|\mathcal{F}| \le \binom{n}{3}+n.$$
\end{theorem}
\begin{proof}
In this case, we have $\alpha_1+\alpha_2+\alpha_3=3n-2(d_1+d_2+d_3) \ge 0$.
Suppose $|\mathcal{F}| > \binom{n}{3}+n+1$. Then there exists a non-zero vector $u \in \mathbb{R}^{\mathcal{F}}$ orthogonal to all $v_S$ for $|S| \in \{0,1,3\}$. We use the previous matrix identity \eqref{eq_key} and compute $u^t M u$ in two different ways. This gives
$$8d_1d_2d_3 \|u\|_2^2=-2(\alpha_1+\alpha_2+\alpha_3)\sum_{|S|=2}\langle u, v_S \rangle^2 \le 0,$$
contradicting that $u$ is a non-zero vector. Therefore we have $|\mathcal{F}| \le \binom{n}{3}+n+1$ for all $n$. 

Next let us show that for $n \ge 24$, it is not possible to have $|\mathcal{F}|=\binom{n}{3}+n+1$. Suppose the contrary. Then from the previous argument, we know that $\{v_S\}$ for $|S|\in \{0,1,3\}$ must form a basis of $\mathbb{R}^{\mathcal{F}}$. In particular, we can find a vector $u \in \mathbb{R}^{\mathcal{F}}$ such that $u$ is orthogonal to $v_S$ whenever $|S| \in \{1,3\}$, and $\langle u, v_{\emptyset}\rangle=1$. In the matrix identity \eqref{eq_key}, we multiply the identity on the right by $u$ and get
\begin{equation}\label{eq_form0}
8d_1d_2d_3 u = c_0 v_{\emptyset} +c_2 \sum_{|S|=2} \langle v_S, u\rangle v_S.
\end{equation}
We let $u^i$ be the vector obtained from $u$ by setting $(u^i)_{x}=u_{x} x_i$. Multiplying the above identity coordinate-wise by $x_i$ gives
\begin{equation} \label{eq_form1}
    8d_1d_2d_3 u^i = c_0 v_{\{i\}} + c_2 \sum_{|S|=2} \langle v_S, u\rangle v_{S \Delta \{i\}}.
\end{equation}
On the other hand, multiplying the matrix identity \eqref{eq_key} with $u^i$ and using $\langle v_S, u^i\rangle = \langle v_{S \Delta \{i\}}, u\rangle$ directly gives
\begin{equation}\label{eq_form2}
8d_1 d_2 d_3 u^i = c_1 \sum_{|S|=1} \langle v_{S \Delta \{i\}}, u \rangle v_S + c_3\sum_{|S|=3} \langle v_{S \Delta \{i\}}, u \rangle v_S.
\end{equation}
Here the $|S| \in \{0,2\}$ terms vanish because after taking symmetric difference with $\{i\}$, they become sets of size $1$ or $3$. 

Now using that $\{v_S\}$ for $|S| \in \{0,1,3\}$ form a basis of $\mathbb{R}^{\mathcal{F}}$ and comparing the coefficients in these two representations \eqref{eq_form1} and \eqref{eq_form2}, we immediately get $c_0=c_1$, and for $i \not\in \{j,k\}$, $c_2 \langle v_{\{j,k\}}, u\rangle = c_3 \langle v_{\{j,k\}}, u\rangle$. Note that $c_2=-2(\alpha_1+\alpha_2+\alpha_3) \le 0<6=c_3$. The latter implies that $u$ is orthogonal to all $\{v_S\}$ for $|S|=2$. Thus by \eqref{eq_form0}, we have
$8d_1d_2d_3 u = c_0 v_{\emptyset}$, and hence $u$ is a vector parallel to the all-one vector. Recall that $\langle u, v_{\emptyset}\rangle=1$. This gives $u=\frac{1}{|\mathcal{F}|} {\bf 1}=\frac{1}{|\mathcal{F}|}v_{\emptyset}.$ Therefore $8d_1d_2d_3=c_0 |\mathcal{F}|$. So clearly we have $c_0>0$ in this case.

Now consider the diagonal entries of the matrix $M$, which gives
\begin{equation}\label{eq_diag}
8d_1d_2d_3=c_0+c_1n+c_2\binom{n}{2}+c_3\binom{n}{3}.
\end{equation}
Plugging in $8d_1d_2d_3=c_0|\mathcal{F}|$, $|\mathcal{F}|=1+n+\binom{n}{3}$ and $c_1=c_0$, we get
$$c_0 \binom{n}{3}=c_2\binom{n}{2}+c_3\binom{n}{3}.$$
Using $c_2=-2(\alpha_1+\alpha_2+\alpha_3)$ and $c_3=6$, this simplifies to
$$(6-c_0)(n-2)=6(\alpha_1+\alpha_2+\alpha_3) \ge 0.$$
Since $c_0>0$ and $c_0$ is an integer by definition, we have $c_0 \in \{1,2, \cdots, 6\}$. Furthermore, it is not hard to see that $\alpha_i=n-2d_i$ has the same parity as $n$. Since
$c_0=-\alpha_1\alpha_2\alpha_3-n(\alpha_1+\alpha_2+\alpha_3)$, we can see that $c_0$ is always even. So $c_0 \in \{2,4,6\}$. If $c_0=2$, then $4(n-2)=6(\alpha_1+\alpha_2+\alpha_3)$, implying $\alpha_1+\alpha_2+\alpha_3$ is even. If $c_0=6$, we have $\alpha_1+\alpha_2+\alpha_3=0$, which is also even. In either case $n$ is even and all $\alpha_i$ thus are also even. But then $c_0=-\alpha_1\alpha_2\alpha_3-n(\alpha_1+\alpha_2+\alpha_3)$ becomes a multiple of $4$, contradicting $c_0 \in \{2, 6\}$.

The only remaining case to check is when $c_0=4$. Then we have $n-2=3(\alpha_1+\alpha_2+\alpha_3)=9n-6(d_1+d_2+d_3)$. This gives
$4n+1=3(d_1+d_2+d_3)$, which violates the assumption that all $d_i$'s are even.
\end{proof}

Our next theorem handles the case when the sum of $d_i$ is ``large''.

\begin{theorem}\label{thm_n_small}
For $n \ge 24$, if $\mathcal{F} \subset \{-1,1\}^n$ is a $3$-distance set with distances from $\{d_1, d_2, d_3\}$ and $d_1+d_2+d_3 > \frac{3n+6}{2}$, then 
$$|\mathcal{F}| < \binom{n}{3}+n.$$
\end{theorem}
\begin{proof}
We define 
$$K(x, y)=\left(\frac{2(d_1+d_2+d_3)}{3}-n-x_1y_1-x_2y_2 \right)P(\langle x, y \rangle).$$
For $x, y \in \mathcal{F}$, when $x \neq y$ we have $P(\langle x, y \rangle)=0$ and thus also $K(x, y)=0$. When $x=y$, we have $P(\langle x, x \rangle)=P(n)=8d_1d_2d_3$, hence
$$K(x, x)=\left(\frac{2(d_1+d_2+d_3)}{3}-n-2 \right) \cdot 8d_1d_2d_3>0.$$
So $K$ satisfies the assumptions in  Lemma \ref{lem_inertia}. From the previous proof, we have 
$$P(\langle x, y \rangle)=c_0 \chi_{\emptyset}(x) \chi_{\emptyset}(y) + c_1 \sum_{|S|=1} \chi_{S}(x)\chi_S(y) + c_2 \sum_{|S|=2} \chi_{S}(x)\chi_S(y)+ c_3 \sum_{|S|=3} \chi_{S}(x)\chi_S(y).$$
Multiplying it by $\frac{2(d_1+d_2+d_3)}{3}-n-x_1y_1-x_2y_2$, we know that $K$ can be written in the form of $\sum_T k_T \chi_T(x) \chi_T(y)$. And $k_T=0$ for $|T| \ge 5$; $k_T=-c_3|T\cap\{1,2\}|\le 0$ for $|T|=4$. For those $T$ with $|T|=3$ and $T \cap \{1,2 \} \neq \emptyset$, we also have
\begin{align*}
k_T&=\left(\frac{2(d_1+d_2+d_3)}{3}-n \right)c_3-c_2|T \cap \{1,2\}|\\
&=\left(4(d_1+d_2+d_3)-6n \right)(1-|T \cap \{1,2\}|) \le 0
\end{align*}
Therefore a strictly positive coefficient $k_T$ can only occur for $|T| \le 2$, or those $T$ with $|T|=3$ disjoint from $\{1, 2\}$. Lemma \ref{lem_inertia} gives
$$|\mathcal{F}| \le 1+\binom{n}{1}+\binom{n}{2}+\binom{n-2}{3}<n+\binom{n}{3}, $$
completing the proof.
\end{proof}
\medskip
\noindent\textbf{Remark.~}
The argument in Theorem~\ref{thm_n_small} also extends to
$s$-distance sets. Suppose that $2\le s\le n$ and that
$\mathcal{F}\subseteq\{-1,1\}^n$ has distances from
$\{d_1,\ldots,d_s\}$, where
$0<d_1<\cdots<d_s\le n$. If $\sum_{i=1}^sd_i>\frac{s(n+2)}{2}$, then
\[
|\mathcal{F}|
\le
\sum_{i=0}^{s-1}\binom{n}{i}
+\binom{n-2}{s}.
\]
In particular, when $n\ge 4s$, this upper bound is strictly less
than
\[
\sum_{\substack{0\le i\le s\\ i\equiv s\pmod* 2}}
\binom{n}{i}.
\]

To see this, let
\[
\alpha_i=n-2d_i,\qquad
P(t)=\prod_{i=1}^s(t-\alpha_i),
\qquad
a=-\frac{1}{s}\sum_{i=1}^s\alpha_i.
\]
The assumption on the distances gives $a>2$. We define
\[
K(x,y)=(a-x_1y_1-x_2y_2)P(\langle x,y\rangle).
\]
As in the proof of Theorem~\ref{thm_n_small}, the restriction of
$K$ to $\mathcal{F}$ is a positive multiple of the identity.
Writing
\[
P(\langle x,y\rangle)
=
\sum_{i=0}^s c_i
\sum_{\substack{T\subseteq[n]\\ |T|=i}}
\chi_T(x)\chi_T(y),
\]
we have $c_s=s!$ and $c_{s-1}=s!a$. Therefore, in the expansion
of $K$, the coefficient of $\chi_T(x)\chi_T(y)$ is
\[
s!a(1-|T\cap\{1,2\}|)
\quad\text{when } |T|=s,
\]
and
\[
-s!|T\cap\{1,2\}|
\quad\text{when } |T|=s+1.
\]
All coefficients with $|T|\ge s+2$ are zero. Hence a strictly
positive coefficient can only occur when $|T|\le s-1$, or when
$|T|=s$ and $T\cap\{1,2\}=\varnothing$. The inertia bound gives the claimed upper bound. Finally, suppose that $n\ge 4s$. By Pascal's identity,
\[
\sum_{\substack{0\le i\le s\\ i\equiv s\pmod* 2}}
\binom{n}{i}
=
\sum_{i=0}^{s}\binom{n-1}{i}.
\]
Consequently,
\[
\begin{aligned}
\sum_{\substack{0\le i\le s\\ i\equiv s\pmod* 2}}
\binom{n}{i}
-\sum_{i=0}^{s-1}\binom{n}{i}
-\binom{n-2}{s}&=
\binom{n-1}{s}
-\sum_{i=0}^{s-2}\binom{n-1}{i}
-\binom{n-2}{s}\\
&=
\binom{n-2}{s-1}
-\sum_{i=0}^{s-2}\binom{n-1}{i}.
\end{aligned}
\]
For $1\le i\le s-1$, we have $\frac{\binom{n-1}{i-1}}{\binom{n-1}{i}}
=\frac{i}{n-i}
\le \frac{s-1}{n-s+1}
<\frac13$. Therefore,
\[
\sum_{i=0}^{s-2}\binom{n-1}{i}
<
\left(\sum_{j=1}^{\infty}3^{-j}\right)
\binom{n-1}{s-1}
=
\frac12\binom{n-1}{s-1}
<
\frac{n-s}{n-1}\binom{n-1}{s-1}
=
\binom{n-2}{s-1}.
\]
Thus the difference is strictly positive, proving the claim.

From now on, we return to the notation for the three-distance case introduced above. Theorems \ref{thm_n_large} and \ref{thm_n_small} take care of the cases when $d_1+d_2+d_3$ is either at most $\frac{3n}{2}$ or strictly greater than $\frac{3n+6}{2}$. So it remains to consider the exceptional interval $d_1+d_2+d_3 \in [\frac{3n+1}{2}, \frac{3n+6}{2}]$. Note that for this interval, $c_2=-2(\alpha_1+\alpha_2+\alpha_3)=4(d_1+d_2+d_3)-6n>0$. In the next three sections we will assume that $d_1+d_2+d_3$ falls into this interval, and have separate discussions based on the parity of $n$.

\section{Within the exceptional interval: the even $n$ case} \label{sec_even}
Before proceeding, we define the trace of a square matrix $A=(A_{ij})_{i,j=1}^m$ by
$\Tr A=\sum_{i=1}^m A_{ii}$. When $n$ is even, the possible values for $d_1+d_2+d_3$ in the exceptional interval are $\frac{3n+2}{2}$, $\frac{3n+4}{2}$, or $\frac{3n+6}{2}$. The following theorem covers all these cases.

\begin{theorem}\label{thm_exceptional_even}
For even $n \ge 24$ and even distances $0<d_1 <d_2<d_3 \le n$ satisfying $d_1+d_2+d_3 \in \{\frac{3n+2}{2}, \frac{3n+4}{2}, \frac{3n+6}{2}\}$, $\mathcal{F} \subset \{-1,1\}^n$ with Hamming distances from $\{d_1, d_2, d_3\}$ satisfies
$$|\mathcal{F}| \le \binom{n}{3}+n.$$
\end{theorem}
\begin{proof}
For some constant $\tau \in (0, 1)$ to be determined later, and a fixed coordinate $i \in [n]$, we let
$$K(x, y)=(1-\tau x_i y_i)P(\langle x, y \rangle).$$
Using the previous identity $$P(\langle x, y \rangle)=c_0 \chi_{\emptyset}(x) \chi_{\emptyset}(y) + c_1 \sum_{|S|=1} \chi_{S}(x)\chi_S(y) + c_2 \sum_{|S|=2} \chi_{S}(x)\chi_S(y)+ c_3 \sum_{|S|=3} \chi_{S}(x)\chi_S(y),$$
we can write $K(x, y)=\sum_T k_T \chi_T(x) \chi_T(y)$, and for $T \subset [n]\setminus \{i\}$, and $|T| \le 3$, 
$$k_T=c_{|T|}-\tau c_{|T|+1}, \qquad k_{T\cup \{i\}}=c_{|T|+1}-\tau c_{|T|}.$$ 
For convenience, we set $c_j=0$ for every $j\ge 4$. 

For each $k$, if $c_k \neq c_{k+1}$, then setting $\tau=1$ would make one of $\{k_T, k_{T\cup \{i\}}\}$ strictly positive and the other strictly negative. Thus the same holds when $\tau \in (0,1)$ is sufficiently close to $1$. On the other hand, if $c_k=c_{k+1} \le 0$, both $k_T$ and $k_{T\cup \{i\}}$ are non-positive for $\tau \in (0,1)$ sufficiently close to $1$. Therefore, if there are no two consecutive $c_k$'s that are equal and strictly positive, we know every pair $\{k_T, k_{T\cup \{i\}}\}$ contributes at most one strictly positive $k_T$. Lemma \ref{lem_inertia} gives 
$$|\mathcal{F}| \le \binom{n-1}{0}+\binom{n-1}{1}+\binom{n-1}{2}+\binom{n-1}{3}=\binom{n}{3}+n.$$
We know $c_2=4(d_1+d_2+d_3)-6n \in \{4, 8, 12\}$ and $c_3=6$, so $c_2$ and $c_3$ are never equal. Likewise $c_3\neq 0=c_4$. So the uncovered cases are: 
$$\textup{(i)}~~c_0=c_1>0; \qquad  \textup{(ii)}~~c_1=c_2>0.$$ 
When $d_1+d_2+d_3=\frac{3n+4}{2}$, we have $c_2=8$. However, all $d_i$'s are even, so $n \equiv 0 \pmod 4$, and thus $\alpha_i=n-2d_i \equiv 0 \pmod 4$. Therefore $c_0=-\alpha_1\alpha_2\alpha_3-n(\alpha_1+\alpha_2+\alpha_3) \equiv 0 \pmod{16}.$ And $c_1=3n-2+(\alpha_1\alpha_2+\alpha_1\alpha_3+\alpha_2\alpha_3)\equiv 2 \pmod 4$. So we can never have $c_0=c_1$ or $c_1=c_2$.

For the remainder of this section, we may assume that $d_1+d_2+d_3 \in \{\frac{3n+2}{2}, \frac{3n+6}{2}\}$, and deal with the two cases when $c_0=c_1>0$ or $c_1=c_2>0$ separately:\medskip\\
\noindent {\bf Case 1: } $c_0=c_1>0$.\smallskip\\
We would like to first show that indeed $c_0=c_1 \ge 12$:\\
\noindent (i) If $d_1+d_2+d_3=\frac{3n+2}{2}$, then $\alpha_1+\alpha_2+\alpha_3=-2$. Comparing the two equalities $c_0=-\alpha_1\alpha_2\alpha_3-n(\alpha_1+\alpha_2+\alpha_3)$ and $c_1=3n-2+(\alpha_1\alpha_2+\alpha_1\alpha_3+\alpha_2\alpha_3)$ gives
$$(\alpha_1+1)(\alpha_2+1)(\alpha_3+1)=1-n<0.$$
Meanwhile we also have $(\alpha_1+1)+(\alpha_2+1)+(\alpha_3+1)=1$. Therefore exactly one of $\alpha_i+1$ is strictly negative and the other two are strictly positive. Since $\alpha_1>\alpha_2>\alpha_3$, we have
$\alpha_3+1<0$ and $\alpha_1+1,\alpha_2+1>0$.
They are all integers, so $\alpha_1,\alpha_2\ge 0$ and
$\alpha_3<0$. Thus $\alpha_1\alpha_2\alpha_3\le 0$. This gives 
$$c_0=-\alpha_1\alpha_2\alpha_3+2n \ge 2n \ge 48.$$

\noindent (ii) If $d_1+d_2+d_3=\frac{3n+6}{2}$, then $\alpha_1+\alpha_2+\alpha_3=-6$. Again using the expressions for $c_0$ and $c_1$, we get
$$(\alpha_1+2)(\alpha_2+2)(\alpha_3+2)=(6n-c_0)+2(c_1-(3n-2))-24+8=c_0-12.$$
Since $3n+6=2(d_1+d_2+d_3) \equiv 0 \pmod 4$, $n \equiv 2 \pmod 4$ and thus $\alpha_i=n-2d_i \equiv 2 \pmod 4.$ Therefore $c_0-12=(\alpha_1+2)(\alpha_2+2)(\alpha_3+2) \equiv 0 \pmod{64}$. Since $c_0$ is assumed to be strictly positive, we have $c_0 \ge 12$.

Therefore, in both cases we have $c_0=c_1 \ge 12$. Using
$$8d_1d_2d_3=c_0+c_1 n + c_2 \binom{n}{2}+c_3 \binom{n}{3},~~~~c_2=-2(\alpha_1+\alpha_2+\alpha_3) \in \{4, 12\},~~c_3=6,$$
the ratio bound from Lemma \ref{lem_ratio} gives
$$|\mathcal{F}| \le \frac{c_0+c_1 n + c_2 \binom{n}{2}+c_3 \binom{n}{3}}{c_0}\le 1+n+\binom{n}{2}+\frac{1}{2}\binom{n}{3} <\binom{n}{3}+n$$
for $n \ge 24$.\\

\noindent{\bf Case 2:} $c_1=c_2>0$. \medskip\\
The $d_1+d_2+d_3 \in \{\frac{3n+2}{2}, \frac{3n+6}{2}\}$ cases correspond to $c_1=c_2=4$ and $c_1=c_2=12$ respectively. We discuss these two subcases below.\\

\noindent (i) $c_1=c_2=4$. We take 
$$K(x, y)=(9-3x_1y_1-3x_2y_2-2x_1x_2y_1y_2) P(\langle x, y\rangle).$$
Obviously $K(x, x)=(9-3-3-2)P(\langle x, x\rangle)=8d_1d_2d_3$. For $x, y \in \mathcal{F}$ and $x \neq y$, $K(x,y)=0$. Again let us write $K(x,y)$ as $\sum_T k_T \chi_T(x) \chi_T(y)$. We have
$$k_T=9c_{|T|}-3c_{|T \Delta \{1\}|}-3c_{|T \Delta \{2\}|}-2c_{|T \Delta \{1,2\}|}.$$
Substituting $c_1=4, c_2=4, c_3=6$, we have that for $i, j \not\in \{1,2\}$,
$$k_{\{i\}}=9c_1-6c_2-2c_3=0,~~k_{\{i,j\}}=9c_2-6c_3=0,$$
$$k_{\{1,i\}}=k_{\{2,i\}}=9c_2-3c_1-3c_3-2c_2=-2.$$
Also when $|T| \ge 4$, clearly $k_T \le 0$. So the strictly positive $k_T$'s can only come from $|T|=3$, or $T \in \{\emptyset, \{1\}, \{2\}, \{1,2\}\}$. Therefore using Lemma \ref{lem_inertia},
$$|\mathcal{F}| \le \binom{n}{3}+4 < \binom{n}{3}+n.$$

\noindent (ii) $c_1=c_2=12$. For this last case, we use a different linear algebraic proof by analyzing the trace of matrices. We let $A_i=\sum_{|S|=i} v_S v_S^T$; then we have the matrix identity
\begin{equation}\label{matrix_1}
8d_1d_2d_3 I = c_0 A_0 + c_1A_1+c_2A_2+c_3A_3=c_0 J+12A_1+12A_2+6A_3,
\end{equation}
where $J$ is the all-one matrix. Note that all these matrices are of dimension $|\mathcal{F}| \times |\mathcal{F}|$. The assertion is immediate if $\mathcal F=\varnothing$, so we may assume that $|\mathcal F|>0$. Now define a new matrix $Q=I-\frac{J}{|\mathcal{F}|}$. This is an orthogonal projection so $Q^2=Q$. Multiplying \eqref{matrix_1} on both sides by $Q$, we have
$$8d_1d_2d_3 Q = 12QA_1Q+12QA_2Q+6QA_3Q.$$
Since every matrix on the right hand side is positive semi-definite, the eigenvalues of $QA_1Q$ are at most $\frac{8d_1d_2d_3}{12}=\frac{2d_1d_2d_3}{3}$. Its rank is at most the number of singletons which is $n$. Therefore
$$\Tr[QA_1Q] \le n \cdot \frac{2d_1d_2d_3}{3}.$$
On the other hand, we have $(A_1)_{x,y}=\langle x, y\rangle$ by definition. So $\Tr[A_1]=\sum_{x \in \mathcal{F}} \langle x, x\rangle = n|\mathcal{F}|$. 
$$\Tr[A_1 J]=\sum_{x, y \in \mathcal{F}} \langle x,y \rangle = n|\mathcal{F}|^2-2\sum_{x, y \in \mathcal{F}}d(x,y).$$
Therefore, we have a second way to compute $\Tr[QA_1Q]$, which gives
\begin{align*}
\Tr[QA_1Q]&=\Tr[A_1Q^2]=\Tr[A_1 Q]=\Tr\left[A_1-\frac{A_1 J}{|\mathcal{F}|}\right]\\
&=n|\mathcal{F}|-n|\mathcal{F}| + \frac{2\sum_{x,y\in \mathcal{F}}d(x,y)}{|\mathcal{F}|}  \ge 2(|\mathcal{F}|-1) d_1. 
\end{align*}
Combining these two results, we get
$$|\mathcal{F}| \le \frac{nd_2d_3}{3}+1.$$
From $c_1=c_2=12$, we get 
$$\alpha_1+\alpha_2+\alpha_3=-6,\qquad \alpha_1\alpha_2+\alpha_1\alpha_3+\alpha_2\alpha_3=14-3n.$$
Together they give 
$$(\alpha_1+2)+(\alpha_2+2)+(\alpha_3+2)=0,
\qquad (\alpha_1+2)^2+(\alpha_2+2)^2+(\alpha_3+2)^2=6n-4.$$
The Cauchy--Schwarz inequality gives 
$$(\alpha_2+2)^2+(\alpha_3+2)^2 \ge \frac{1}{2}((\alpha_2+2)+(\alpha_3+2))^2=\frac{1}{2}(\alpha_1+2)^2.$$
Thus $(\alpha_1+2)^2 < 4n$ and $\alpha_1<2\sqrt{n}-2$. This immediately gives
$$d_2+d_3=\frac{n-\alpha_2}{2}+\frac{n-\alpha_3}{2}=n+3+\frac{\alpha_1}{2}<n+\sqrt{n}+2.$$
Therefore 
$$|\mathcal{F}| \le 1+\frac{nd_2d_3}{3} \le 1+\frac{n(n+\sqrt{n}+2)^2}{12},$$
which is strictly less than $n+\binom{n}{3}$ for $n \ge 24.$
\end{proof}

\section{Within the exceptional interval: the odd $n$ case}\label{sec_odd}
When $n$ is odd, the remaining possible distance sums are $d_1+d_2+d_3 \in \{\frac{3n+1}{2}, \frac{3n+3}{2}, \frac{3n+5}{2}\}$. We first deal with the first and last cases, i.e. $d_1+d_2+d_3$ is equal to $\frac{3n+1}{2}$ or $\frac{3n+5}{2}$. 
\begin{theorem}
Let $n\ge 24$ be odd, and let
$0<d_1<d_2<d_3\le n$ be even integers satisfying $d_1+d_2+d_3
\in
\left\{\frac{3n+1}{2},\frac{3n+5}{2}\right\}.$
Suppose that $\mathcal F\subseteq\{-1,1\}^{n}$ has all pairwise
distances in $\{d_1,d_2,d_3\}$. Then
\[
|\mathcal F|\le n+\binom{n}{3}.
\]
\end{theorem}
\begin{proof}
In these two cases we have $n \equiv 1 \pmod 4$, $c_2 \in \{2, 10\}$, $c_3=6$. Using the same function 
$$K(x,y)=(1-\tau x_i y_i)P(\langle x, y \rangle)$$ and the same arguments as in the proof for the even $n$ case in Section \ref{sec_even}, we know that if there are no consecutive $c_k=c_{k+1}>0$, then $|\mathcal{F}| \le \binom{n}{3}+n$. Clearly $c_2 \neq c_3$. We also know that $\alpha_i=n-2d_i \equiv 1 \pmod 4$, so $c_1\equiv 3n-2+3 \equiv 0 \pmod 4$ and thus $c_1 \neq c_2$. The only case to check is $c_0=c_1>0$. 

When $d_1+d_2+d_3=\frac{3n+1}{2}$, calculations show that 
$$\prod_{i=1}^3 (\alpha_i+1)=-2(n-1)<0,\qquad \sum_{i=1}^3 (\alpha_i+1)=2>0.$$ So one of $\alpha_i+1$ is strictly negative and the other two are strictly positive. Therefore $\alpha_1\alpha_2\alpha_3 \le 0$ and thus $c_0=-\alpha_1\alpha_2\alpha_3-n(\alpha_1+\alpha_2+\alpha_3) \ge n$. 

When $d_1+d_2+d_3 =\frac{3n+5}{2}$, we have $\prod_{i=1}^3 (\alpha_i+1)=2(n-1)>0$, $\sum_{i=1}^3 (\alpha_i+1)=-2<0$. So two of $\alpha_i+1$ are strictly negative and one is strictly positive. Since $\alpha_i+1 \equiv 2 \pmod 4$, the negative $\alpha_i+1$ are at most $-2$. Therefore $\prod_{i=1}^3 (\alpha_i+3) \ge 0$. Expanding this gives $c_0 \ge 2n+6$. Therefore Lemma \ref{lem_ratio} gives
\begin{align*}
|\mathcal{F}| \le \frac{c_0+c_1n+c_2\binom{n}{2}+c_3\binom{n}{3}}{c_0} \le 1+n+\frac{10}{n}\binom{n}{2}+\frac{6}{n}\binom{n}{3}=n^2+3n-2 < \binom{n}{3}+n
\end{align*}
for $n \ge 24$. 
\end{proof}

Therefore it remains to check the middle case $d_1+d_2+d_3=\frac{3n+3}{2}$. In this case we have 
$$c_0=3n-\alpha_1\alpha_2\alpha_3,\qquad c_1=3n-2+\sum_{1\le i <j \le 3} \alpha_i \alpha_j, \qquad c_2=c_3=6.$$
We introduce one additional parameter 
$$\rho=c_1-c_0=\alpha_1\alpha_2\alpha_3+\left(\sum_{1 \le i<j \le 3}\alpha_i\alpha_j\right)-2 = \prod_{i=1}^3(\alpha_i+1).$$
We make some quick number-theoretic observations here. First we have  $\alpha_i \equiv n \pmod 4$, so $-3=\sum_{i=1}^3 \alpha_i \equiv 3n \pmod 4$, so $n \equiv 3 \pmod 4$. This gives $n \ge 27$ and $\alpha_i \equiv 3 \pmod 4$. Therefore $\rho=\prod_{i=1}^3 (\alpha_i+1) \equiv 0 \pmod{64}$. Also $c_0 \equiv 2 \pmod 4$, $c_1\equiv 2 \pmod 4$.

The next theorem covers the middle case when $c_1 \ge 6$. The techniques used in this proof do not apply to the $c_1<6$ range. We will address that final and most challenging case in the next section.

\begin{theorem}
Let $n\ge 24$ be odd, and let
$0<d_1<d_2<d_3\le n$ be even integers satisfying $d_1+d_2+d_3=\frac{3n+3}{2}$. Suppose that $\mathcal F\subseteq\{-1,1\}^{n}$ has all pairwise
distances in $\{d_1,d_2,d_3\}$. If $c_1\ge 6$, then
\[
|\mathcal F|\le n+\binom{n}{3}.
\]
\end{theorem}
\begin{proof}
We discuss the following four cases.\\

\noindent \textbf{Case 1:} $c_1 \ge 6$ and $\rho \le 0$, except $c_0=c_1=6$.\smallskip

Since $\rho=c_1-c_0$, we have $c_0 \ge c_1 \ge 6$. Since we have excluded the $c_0=c_1=6$ case, we have $c_0>6$ and thus since $c_0 \equiv 2 \pmod 4$, we have $c_0 \ge 10$. The ratio bound from Lemma \ref{lem_ratio} gives
$$|\mathcal{F}| \le 1+\frac{c_1}{c_0} n + \frac{6}{c_0}\left(\binom{n}{2}+\binom{n}{3}\right) \le 1+n+\frac{3}{5}\left(\binom{n}{2}+\binom{n}{3}\right)<\binom{n}{3}+n$$
for $n \ge 24$.\\

\noindent \textbf{Case 2:} $c_1 > 6$ and $\rho>0$. \smallskip

Consider the vectors $\{v_S\}$ previously defined, with $v_S \in \mathbb{R}^{\mathcal{F}}$  and $(v_S)_x=\prod_{i \in S} x_i$. If $\{v_S\}_{S: |S| \in \{1, 3\}}$ span $\mathbb{R}^{\mathcal{F}}$, then the dimension argument immediately gives our desired inequality $|\mathcal{F}| \le \binom{n}{3}+n$. So we may assume a non-zero vector $u$ is orthogonal to $v_S$ when $|S| \in \{1, 3\}$. Recall that in the proof of Theorem \ref{thm_exceptional_even}, we introduced the matrix $A_i$ with $(A_i)_{x,y}=\sum_{S:\,|S|=i}\chi_S(x)\chi_S(y).$
Therefore, for every $x,y\in\mathcal F$,
\[
\langle x,y\rangle (A_i)_{x,y}
=
(i+1)(A_{i+1})_{x,y}
+
(n-i+1)(A_{i-1})_{x,y}.
\]
Define the matrix $K=(K(x,y))_{x,y\in\mathcal F}$ by
\[
K(x,y)
=
(\langle x,y\rangle-n)P(\langle x,y\rangle).
\]
The matrix $K$ is zero: if $x\ne y$, then
$P(\langle x,y\rangle)=0$, while if $x=y$, then
$\langle x,y\rangle-n=0$. On the other hand, the preceding entrywise identity gives
\begin{align*}
K={}&
c_0A_1+c_1(2A_2+nA_0)
+c_2(3A_3+(n-1)A_1)+c_3(4A_4+(n-2)A_2)
-n(c_0A_0+c_1A_1+c_2A_2+c_3A_3)\\
={}&
n\rho A_0
-\bigl(\rho+(n-1)(c_1-6)\bigr)A_1
+2(c_1-6)A_2-6(n-3)A_3+24A_4.
\end{align*}
Therefore the right hand side is an all-zero matrix on $\mathcal{F} \times \mathcal{F}$. Treating it as a quadratic form and plugging in $u$, we have
$$0=n \rho \langle u, v_{\emptyset}\rangle^2+ 2(c_1-6) \sum_{|S|=2} \langle u, v_S\rangle^2 + 24 \sum_{|S|=4} \langle u, v_S\rangle^2$$
Since all the terms involved are non-negative and $\rho>0$, $c_1-6>0$, we know $\langle u, v_{S}\rangle=0$ for $|S| \in \{0,2\}$. Therefore multiplying the equality $8d_1d_2d_3 I = c_0A_0+c_1A_1+c_2A_2+c_3A_3$ by $u$ on its right, we have
$u=0$, contradicting the assumption that $u$ is a non-zero vector.\\

\noindent \textbf{Case 3:} $c_1=c_0=6$. \smallskip

In this case we have $\rho=c_1-c_0=0$ while $\rho=\prod_{i=1}^3 (\alpha_i+1)$. Since $\sum_{i=1}^3(\alpha_i+1)=0$, we can assume
$\alpha_1+1=\alpha$, $\alpha_2+1=0$, and
$\alpha_3+1=-\alpha$ for some integer $\alpha>0$. Using $\alpha_i=n-2d_i$, we have
$$d_1=\frac{n+1-\alpha}{2},\qquad d_2=\frac{n+1}{2},\qquad d_3=\frac{n+1+\alpha}{2},\qquad \alpha^2=1+\alpha_1\alpha_2\alpha_3=3n-5.$$
Now consider the $|\mathcal{F}| \times |\mathcal{F}|$ matrix $K$ with 
$$K(x, y)=\frac{(d_1-d(x,y))(d_3-d(x,y))}{(d_1-d_2)(d_3-d_2)}.$$
Note that the columns of $K$ are linear combinations of evaluations of $n$-variable multilinear polynomials of degree at most $2$, so the rank of $K$ is at most $1+n+\binom{n}{2}$. On the other hand, 
$$K= \left(1-\frac{(n+1)^2}{\alpha^2}\right)I+ B,$$ where $B$ is some $0/1$ binary symmetric matrix for distance $d_2.$ Therefore either $|\mathcal{F}| \le 1+n+\binom{n}{2}<\binom{n}{3}+n$, or $0$ is an eigenvalue of $K$ and thus $-1+\frac{(n+1)^2}{\alpha^2}$ is an eigenvalue of $B$. Since all entries of $B$ are integers, its characteristic polynomial is an integer-coefficient polynomial with leading coefficient $1$. Hence $-1+\frac{(n+1)^2}{\alpha^2}$ must also be an integer. This gives $\alpha \mid n+1$. However $\alpha^2=3n-5=3(n+1)-8$. So $\alpha \mid 8$. Recall that $\alpha=\alpha_1+1 \equiv 0 \pmod 4$. The only possibilities are $\alpha=4$ or $\alpha=8$, which correspond to $n=7$ and $n=23$ respectively, contradicting $n \ge 24$.\\

\noindent \textbf{Case 4:} $c_1=6$ and $\rho>0$ (since $\rho \equiv 0 \pmod{64}$ we have $\rho \ge 64$). \smallskip

We take 
$$K(x,y)=\left(\langle x, y\rangle^2 + \left(2+\frac{\rho}{3(n-3)}\right) \langle x, y\rangle+ 6-7n-\frac{\rho}{n-3}\right)P(\langle x, y \rangle).$$
The coefficients are chosen so that the coefficients of $A_2$ and $A_3$ vanish. Calculations give 
$$gI=g_0A_0+g_1A_1+g_4A_4+g_5A_5,$$
where 
\begin{align*}
g&=K(x,x)=8d_1d_2d_3 \left(n^2+\left(2+\frac{\rho}{3(n-3)}\right) n+ 6-7n-\frac{\rho}{n-3}\right)\\
&=8d_1d_2d_3\left((n-2)(n-3)+\frac{\rho}{3}\right) = \left(c_0+c_1n+c_2\binom{n}{2}+c_3\binom{n}{3}\right)\left((n-2)(n-3)+\frac{\rho}{3}\right)\\
&=\left((n+1)(n^2-n+6)-\rho\right)\left((n-2)(n-3)+\frac{\rho}{3}\right)
\end{align*}
$$g_0=6(n-2)(n-3)+(6n-4)\rho+\frac{\rho^2}{n-3},$$
$$g_1=6(n-2)(n-3)-\frac{\rho^2}{3(n-3)},\qquad g_4=120+\frac{8\rho}{n-3},\qquad g_5=120.$$
Obviously $g_0, g_4, g_5>0$. We will show that $g_1$ is positive too. Note that
$$\sum_{i=1}^3 (\alpha_i+1)^2=\left(\sum_{i=1}^3 \alpha_i\right)^2-2\sum_{1\le i<j\le 3}\alpha_i\alpha_j + 2\sum_{i=1}^3 \alpha_i +3 =6n-10.$$
Therefore by the AM-GM inequality,
$$\rho^2=\prod_{i=1}^3 (\alpha_i+1)^2 \le \left(\frac{6n-10}{3}\right)^3<8n^3 \le 18(n-2)(n-3)^2$$
for $n \ge 27$. This immediately gives $g_1>0$. Before applying the ratio bound, we prove that $g_0>\frac{6(n+1)}{n-2}((n-2)(n-3)+\frac{\rho}{3})$. Calculation shows that the difference between the left and right hand sides is equal to (using $\rho\ge 64 \ge 4$):
$$-18(n-3)+\frac{6\rho(n^2-3n+1)}{(n-2)}+\frac{\rho^2}{n-3} \ge \frac{6(n^2+3n-14)}{n-2} >0.$$

Therefore the ratio bound in Lemma \ref{lem_ratio} gives 
$$|\mathcal{F}| \le \frac{g}{g_0}<\frac{(n+1)(n^2-n+6)-\rho}{\frac{6(n+1)}{n-2}}<\frac{(n-2)(n^2-n+6)}{6}=\binom{n}{3}+n-2< \binom{n}{3}+n,$$
which completes the proof of the last case.
\end{proof}

\section{The last case: $n$ odd, $d_1+d_2+d_3=\frac{3n+3}{2}$ and $c_1<6$}\label{sec_difficult}

In this section, we will tackle the last and perhaps the most challenging case. Throughout this section, we will assume that $n \ge 24$ is odd, $d_1+d_2+d_3=\frac{3n+3}{2}$, and $c_1<6$. 

Let us first provide a proof sketch to illustrate the key ideas. Let $O$ be the $|\mathcal F|\times\left(n+\binom n3\right)$
matrix whose first $n$ columns are
$v_{\{1\}},\dots,v_{\{n\}}$, followed by the columns $v_T$
for all three-element subsets $T\subseteq[n]$ in any fixed order.
The rows are indexed by $x\in\mathcal F$. Thus the entry in
row $x$ and column $S$ is
\[
O_{x,S}=(v_S)_x=\chi_S(x)=\prod_{i\in S}x_i.
\]

With this column ordering, we regard $O$ as a linear map $O:\mathbb R^n\oplus\mathbb R^{\binom n3}
\longrightarrow\mathbb R^{\mathcal F}$.
Writing a vector in its domain as $(\ell,p)$, where
$\ell=(\ell_i)_{i=1}^n$ and $p=(p_T)_{|T|=3}$, we have $O(\ell,p)
=\sum_{i=1}^n\ell_i v_{\{i\}}
+\sum_{|T|=3}p_Tv_T.$
Equivalently, $\bigl(O(\ell,p)\bigr)_x
=\sum_{i=1}^n\ell_i x_i
+\sum_{|T|=3}p_T\chi_T(x)$ for every $x\in\mathcal F$. We then let $U=\ker O^T$ and $V=\ker O$. Then by comparing the row rank and column rank of $O$, we have
$$\binom{n}{3}+n-\dim V = |\mathcal{F}|-\dim U.$$
Our goal is to show that $\dim U \le \dim V$; then immediately $|\mathcal{F}| \le \binom{n}{3}+n$.

To compare the dimensions of $U$ and $V$, we will construct linear maps $T_i: U \rightarrow V \subset \mathbb{R}^n \oplus \mathbb{R}^{\binom{n}{3}}$ for $i=1, \cdots, n$, such that for a carefully chosen symmetric bilinear form $g$ on $\mathbb{R}^n \oplus \mathbb{R}^{\binom{n}{3}}$, 
$$g(T_i(u), T_j(v))=g(T_j(u), T_i(v)).$$
We will prove a dimension lemma (Lemma \ref{lem_kw}) showing that if there are subspaces $K \subset U$, $W \subset V$, such that $g(w, w')=0$ whenever $w, w' \in W$, and $W$ contains every $T_i(K)$, then $\dim K \le \dim W$. Finally, to compare $\dim U$ and $\dim V$, in
Lemma~\ref{lem_uv} we construct subspaces $K\subseteq U$
and $W\subseteq V$ satisfying the hypotheses of the preceding
lemma, together with a nonnegative integer $r$, such that
\[
\dim U=r+\dim K,
\qquad
\dim V\ge r+\dim W.
\]
Here $r$ is the rank of a linear map constructed in the proof.
The preceding lemma gives $\dim K\le\dim W$, and hence
\[
\dim U=r+\dim K\le r+\dim W\le\dim V.
\]
By the preceding discussion, this implies
$|\mathcal F|\le\binom n3+n$.

In the remainder of this section we establish a series of lemmas needed to complete this proof using the approaches outlined above. First we explain how to construct $T_i: U \rightarrow V$. For each $u \in U=\ker O^T \subset \mathbb{R}^{\mathcal{F}}$, we associate with it a symmetric $n \times n$ matrix $Q_u$, and a $\binom{n}{3} \times n$ matrix $C_u$ whose rows are labeled by $3$-subsets $T$ of $[n]$, by 
$$(Q_u)_{i,i}=(c_0-c_1) \langle u, v_{\emptyset}\rangle,\qquad (Q_u)_{i,j}=(6-c_1)\langle u, v_{\{i,j\}}\rangle,~~~\textup{for~}i \neq j.$$
$$(C_u)_{T,i}=\begin{cases}
\langle u, v_{T \cup \{i\}}\rangle, \qquad&  i\not\in T,\\
0, \qquad& i \in T.
\end{cases}$$

\begin{lemma}\label{lem_inj}
The mapping $u \rightarrow Q_u$ is injective.
\end{lemma}
\begin{proof}
From \eqref{eq_key}, we have
$$8d_1d_2d_3 \|u\|_2^2=\sum_{i=0}^3 c_i \sum_{|S|=i}\langle u, v_S\rangle^2.$$
Since $u \in \ker O^T$, we have $\langle u, v_S \rangle=0$ whenever $|S| \in \{1,3\}$. Therefore
$$8d_1d_2d_3 \|u\|_2^2=c_0 \langle u, v_{\emptyset}\rangle^2+ c_2 \sum_{|S|=2}\langle u, v_S\rangle^2=c_0 \langle u, v_{\emptyset}\rangle^2+ 6 \sum_{|S|=2}\langle u, v_S\rangle^2.$$
Suppose $Q_u=0$. If $c_0 \neq c_1$, then using both $c_0-c_1 \neq 0$ and $6-c_1 \neq 0$, we have $\langle u, v_S\rangle =0$ for $|S| \in \{0, 2\}$. The above equality immediately gives $\|u\|_2^2=0$ and thus $u=0$. If $c_0=c_1$, recall that this gives $\rho=\prod_{i=1}^3(\alpha_i+1)=0$. Together with $\sum_{i=1}^3 (\alpha_i+1)=0$, we have $\alpha_1=\alpha-1, \alpha_2=-1, \alpha_3=-\alpha-1$, and then $c_0=c_1=3n+1-\alpha^2$. Since $c_1 < 6$ and $c_1 \equiv 2 \pmod 4$, we have $c_1 \le 2$. However if $c_1=2$, we have $3n+1-\alpha^2=c_1=2$ and this implies $\alpha^2 \equiv -1 \pmod 3$, which is clearly impossible. Therefore $c_0=c_1 \le -2$. Note that since $6-c_1 \neq 0$, from $Q_u=0$ we still have $\langle u, v_S\rangle=0$ whenever $|S|=2$ . Then the equality above gives 
$$8d_1d_2d_3 \|u\|_2^2 = c_0 \langle u, v_{\emptyset}\rangle^2,$$
where the right hand side is non-positive since $c_0 \le -2$, and the left hand side is non-negative. So the only possibility is $u=0$. So in either case, $Q_u=0$ implies $u=0$, proving the injectivity of $u \rightarrow Q_u$.
\end{proof}

The mappings $T_i: U \rightarrow V$ are now defined as follows:
\[
T_i(u)=(Q_u e_i,-6C_u e_i),
\]
where $e_i$ is the $i$-th standard unit vector. An alternative way is to view $T=(T_i)_{i \in [n]}$ as a mapping from $U$ to the space of linear maps from $\mathbb{R}^n$ to $V$. But we will not adopt this convention below, and try to keep our proof as simple as possible.

In the space $\mathbb{R}^n \oplus \mathbb{R}^{\binom{n}{3}}$ where $V$ lies, we define the following symmetric bilinear form for $\ell, \ell' \in \mathbb{R}^n$, $p, p' \in \mathbb{R}^{\binom{n}{3}}$:
$$g((\ell, p), (\ell', p'))=-\frac{1}{6-c_1} \langle \ell, \ell'\rangle + \frac{1}{6}\langle p, p'\rangle.$$

\begin{lemma}\label{lem_compatibility}
The mappings $T_i$ indeed map $U$ to $V=\ker O$, and the bilinear form $g$ satisfies the following compatibility relation:
$$g(T_i(u), T_j(v))=g(T_j(u), T_i(v)).$$
\end{lemma}
\begin{proof}
First we check that $T_i(u) \in V=\ker O$. 
Multiplying \eqref{eq_key} by $u$ and using that
$u\in U=\ker O^T$ is orthogonal to $v_S$ for
$|S|\in\{1,3\}$, we obtain
\begin{align*}8d_1d_2d_3 u &= c_0 \langle u, v_{\emptyset}\rangle v_{\emptyset} + c_2 \sum_{|S|=2} \langle u, v_S\rangle v_S=c_0 \langle u, v_{\emptyset}\rangle v_{\emptyset} + 6 \sum_{|S|=2} \langle u, v_S\rangle v_S.
\end{align*}
Let $a_S(u)=\langle u, v_S\rangle=\sum_x u_x \chi_S(x)$; then the right hand side becomes $c_0a_{\emptyset}(u)v_{\emptyset}+6\sum_{|S|=2} a_S(u)v_S$. Note that $x_i \chi_S(x)=\chi_{S\Delta \{i\}}(x)$. Therefore, by restricting to the $x$-th coordinate for $x \in \mathcal{F} \subset \{-1, 1\}^n$ and then multiplying both sides by $x_i$, we have
$$8d_1d_2d_3 u_x x_i = c_0 a_{\emptyset}(u) x_i+ 6 \sum_{j \neq i} a_{\{i,j\}}(u)x_j + 6 \sum_{|T|=3, i \in T} a_{T \setminus \{i\}}(u) \chi_T(x).$$
Define from $u$ a new vector $w$ with $w_x=x_i u_x$. We have 
$$a_S(w)=\sum_{x} w_x \chi_S(x)=\sum_x x_i u_x\chi_S(x)=\sum_x u_x \chi_{S \Delta \{i\}}(x).$$
Therefore $a_S(w)=a_{S \Delta \{i\}}(u)$. Applying \eqref{eq_key} to a vector $w$  again, we get
\begin{align*}
8d_1d_2d_3 u_x x_i &= 8d_1d_2d_3 w_x \\
&=  c_1 a_{\emptyset}(u) x_i+c_1 \sum_{j \neq i} a_{\{i,j\}}(u) x_j + 6 \sum_{|T|=3, i \in T} a_{T\setminus \{i\}}(u) \chi_T(x) + 6 \sum_{|T|=3, i \not\in T} a_{T \cup \{i\}}(u) \chi_T(x).
\end{align*}
Subtracting these two identities involving $8d_1d_2d_3 u_x x_i$, after cancelling some common terms on the right hand side, we get
$$0=(c_0-c_1) a_{\emptyset}(u) x_i + (6-c_1) \sum_{j \neq i} a_{\{i,j\}}(u)x_j - 6 \sum_{|T|=3, i \not\in T} a_{T \cup \{i\}}(u) \chi_T(x).$$
This is exactly equivalent to $T_i(u) \in \ker O$. So $T_i$ is a well-defined mapping from $U$ to $V$.

It remains to verify the compatibility relation. For arbitrary $(\ell,p)\in V$, we have
\begin{align*}
g(T_i(u), (\ell, p))&=-\frac{1}{6-c_1} \langle Q_u e_i, \ell\rangle + \frac{1}{6}\langle -6C_u e_i, p\rangle=-\frac{1}{6-c_1} \langle Q_u e_i, \ell\rangle -\langle C_u e_i, p\rangle.
\end{align*}
By the definitions of the matrices $Q_u$ and $C_u$,
$$\langle Q_u e_i, \ell\rangle = (c_0-c_1) \ell_i a_{\emptyset}(u) + (6-c_1)\sum_{j \neq i} \ell_j a_{\{i,j\}}(u).$$
$$\langle C_u e_i, p\rangle = \sum_{|T|=3, i \not\in T}a_{T \cup \{i\}}(u) p_T.$$
Therefore 
\begin{align} \label{eq_g1}
g(T_i(u), (\ell, p))=\frac{c_1-c_0}{6-c_1} \ell_i a_{\emptyset}(u) -\sum_{j \neq i} \ell_j a_{\{i,j\}}(u)-\sum_{|T|=3, i \not\in T}p_T a_{T \cup \{i\}}(u).
\end{align}
Since $(\ell,p)\in V=\ker O$, we have, for every $x\in\mathcal F$,
\[
\sum_{j=1}^{n}\ell_jx_j
+
\sum_{|T|=3}p_T\chi_T(x)
=0.
\]
Multiplying this identity by $u_xx_i$ and summing over
$x\in\mathcal F$, we obtain
\begin{align}\label{eq_g2}
0={}&\ell_i a_{\emptyset}(u)
+\sum_{j\ne i}\ell_j a_{\{i,j\}}(u)
+\sum_{\substack{|T|=3\\ i\in T}}
p_Ta_{T\setminus\{i\}}(u) +
\sum_{\substack{|T|=3\\ i\notin T}}
p_Ta_{T\cup\{i\}}(u).
\end{align}
Adding \eqref{eq_g1} and \eqref{eq_g2}, we obtain
$$g(T_i(u), (\ell, p))=\frac{6-c_0}{6-c_1} \ell_i a_{\emptyset}(u) + \sum_{|T|=3, i\in T}
p_Ta_{T\setminus\{i\}}(u).$$
If $i=j$, the desired compatibility relation is immediate.
Hence we may assume that $i\ne j$.

Let $(\ell,p)=T_j(v)=(Q_ve_j,-6C_ve_j)$. Then
\[
\ell_i=(Q_v)_{ij}=(6-c_1)a_{\{i,j\}}(v),
\]
and
\[
p_T=
\begin{cases}
-6a_{T\cup\{j\}}(v), & j\notin T,\\
0, & j\in T.
\end{cases}
\]
Substituting these expressions into the preceding identity gives
\[
\begin{aligned}
g(T_i(u),T_j(v))
={}&(6-c_0)a_{\emptyset}(u)a_{\{i,j\}}(v)-6\sum_{\{k,\ell\}\subseteq[n]\setminus\{i,j\}}
a_{\{k,\ell\}}(u)\,
a_{\{i,j,k,\ell\}}(v).
\end{aligned}
\]
The right-hand side is unchanged when $i$ and $j$ are interchanged.
Therefore
\[
g(T_i(u),T_j(v))=g(T_j(u),T_i(v)).
\]
\end{proof}

Our next lemma is purely linear algebraic. Since it is not directly related to our extremal combinatorial problem, we defer its proof to Appendix \ref{app_1} for interested readers. For a linear map $A:X\to Y$, we write
$\im A=\{A(x):x\in X\}$ for its image.
For a matrix $A$, its image is its column space.

\begin{lemma}\label{lem_la}
Let $\mathcal{Q}$ be a linear space of symmetric $n \times n$ matrices. Suppose $R$ is a linear map from $\mathcal{Q}$ to the linear space of $n \times n^3$ matrices, i.e. each matrix $Q \in \mathcal{Q}$ is mapped to an $n \times n^3$ matrix $R_Q$. Suppose \smallskip\\
\textup{(i)} the entry $(R_Q)_{i,(j,k,\ell)}$ is invariant under every permutation of $i, j, k, \ell$,\\ 
\textup{(ii)} $R_Q R_P^T=QP$ for $P, Q \in \mathcal{Q}$,\smallskip\\
then  $$\dim \mathcal{Q} \le \dim \left(\sum_{Q \in \mathcal{Q}} \im Q\right).$$
Here $\sum$ means the sum of the column spaces of the matrices $Q \in \mathcal{Q}$.
\end{lemma}

Our next lemma gives a dimension bound using a totally isotropic subspace $W$ of the bilinear form $g$ we defined earlier.

\begin{lemma}\label{lem_kw}
Suppose subspaces $K \subset U$, $W \subset V$ satisfy\smallskip\\
\textup{(i)} $g(w, w')=0$ whenever $w, w' \in W$,\\
\textup{(ii)} $T_i(K) \subset W$ for every $i \in [n]$, \smallskip\\
then 
$$\dim K \le \dim W.$$
\end{lemma}
\begin{proof}
If $K=\{0\}$ then the conclusion trivially holds. Otherwise, assumption (ii) gives $T_i(u)\in W$ for every $u\in K$ and $i\in[n]$. Therefore, by (i), we have $g(T_i(u),T_j(v))=0$ for any $u,v\in K$ and all $i,j\in[n]$.

Using the definitions $T_i(u)=(Q_u e_i, -6C_u e_i)$ and $g((\ell,p),(\ell', p'))=-\frac{1}{6-c_1} \langle \ell, \ell'\rangle + \frac{1}{6}\langle p, p'\rangle$, we have
$$0=-\frac{1}{6-c_1}\langle Q_u e_i, Q_v e_j \rangle +6\langle C_u e_i, C_v e_j \rangle$$
for any $i, j\in [n]$. In particular, this implies as $n \times n$ matrices,
$$6C_u^T C_v=\frac{1}{6-c_1} Q_u^T Q_v.$$
We set $\mathcal{Q}=\{Q_u: u \in K\}$. The injectivity of $u \rightarrow Q_u$ from Lemma \ref{lem_inj} allows us to define $n \times n^3$ matrices $R_{Q_u}$ as
$$(R_{Q_u})_{i,(j,k,\ell)}=\begin{cases}\sqrt{6-c_1}a_{\{i,j,k,\ell\}}(u),\qquad & i,j,k,\ell \textup{~all distinct,}\\
0,\qquad & \textup{otherwise}.
\end{cases}$$
Then clearly $R_{Q_u}$ is invariant under permutations of $\{i,j,k,\ell\}$. Furthermore,
\begin{align*}
(R_{Q_u} R_{Q_v}^T)_{i,i'}&=6(6-c_1)\sum_{j<k<\ell} a_{\{i,j,k,\ell\}}(u)a_{\{i',j,k,\ell\}}(v)=6(6-c_1) \sum_{T: |T|=3} a_{T \cup \{i\}}(u) a_{T\cup \{i'\}}(v)\\
&=6(6-c_1)  \sum_{|T|=3} (C_u)_{T,i} (C_v)_{T,i'} =6(6-c_1) (C_u^T C_v)_{i,i'} = (Q_u^T Q_v)_{i,i'},
\end{align*}
where the terms with $i\in T$ or $i'\in T$ vanish, since
$u,v\in U=\ker O^T$ and $|T|=3$. Therefore $R_{Q_u}R_{Q_v}^T=Q_u^T Q_v$ for any $u,v \in K$ (i.e. $Q_u, Q_v \in \mathcal{Q}$), hence by Lemma \ref{lem_la}, 
$$\dim K = \dim \mathcal{Q} \le \dim \left(\sum_{u \in K} \im Q_u\right) \le \dim W.$$
The first equality $\dim K=\dim \mathcal{Q}$ comes from the injectivity of $u \rightarrow Q_u$. To see the second inequality, we let $P$ be the projection from $V \subset \mathbb{R}^n \oplus \mathbb{R}^{\binom{n}{3}}$ to the first part, i.e. $P((\ell, p))=\ell$. Note that $\im Q_u$ is spanned by $\{Q_u e_i\}$ for $i \in [n]$. But $Q_u e_i=P(T_i(u)) \in P(T_i(K)) \subset P(W)$. Therefore clearly
$$\dim \left(\sum_{u \in K} \im Q_u\right) \le \dim P(W) \le \dim W,$$
completing the proof.
\end{proof}
We state the last lemma needed to complete the proof. It compares the dimension of $U$ and that of $V$ by constructing $K \subset U, W \subset V$ satisfying the conditions in Lemma \ref{lem_kw}.
\begin{lemma}\label{lem_uv}
For the subspaces $U=\ker O^T \subset \mathbb{R}^{\mathcal{F}}$ and $V=\ker O \subset \mathbb{R}^n \oplus \mathbb{R}^{\binom{n}{3}}$ we defined earlier, 
$$\dim U \le \dim V.$$
And therefore,
$$|\mathcal{F}| \le \binom{n}{3}+n.$$
\end{lemma}
\begin{proof}
As we described earlier, our strategy is to construct $K \subset U$ and $W \subset V$, such that $\dim U=r+\dim K$, $\dim V \ge \dim W+r$, and $K, W$ satisfy the conditions in Lemma \ref{lem_kw}, and this would complete the proof.

Let $V^*$ be the space of linear functions from $V$ to $\mathbb{R}$. For each $i \in [n]$, define $B_i: U \rightarrow V^*$ such that 
$$B_i(u)(v)=g(T_i(u), v).$$
We choose coefficients $\{a_i\}_{i \in [n]}$ such that the rank of $B=\sum_{i=1}^n a_i B_i$ is maximized. Suppose the rank of $B$ is equal to $r$. Let $T=\sum_{i=1}^n a_i T_i$, which is also a mapping from $U$ to $V$. Then by the bilinearity of $g$,
$$g(T(u), v)=\sum_{i=1}^n a_i g(T_i(u), v)=\sum_{i=1}^n a_i B_i(u)(v)=B(u)(v).$$
We define $K$ and $W$ (to which we will apply Lemma \ref{lem_kw} later) as follows:
$$K=\ker B \subset U,\qquad W=\spa \{T_i u: u \in K\}_{i\in [n] } \subset V.$$

Since the rank of $B$ is $r$, we know $\dim U=\dim K +r$. Next let us show $\dim W \le \dim V-r$. For every $u \in K, w \in U$, we have $B(w)(T_i(u))=g(T(w), T_i(u))=g(T_i(w), T(u))=B(u)(T_i(w))$. The second last equality comes from the compatibility part of Lemma \ref{lem_compatibility}. Since $u \in K=\ker B$, we have $B(u)=0$, so $B(w)(T_i(u))=0$ too. Therefore every function in $\im B$ vanishes on $T_i(u)$ for every $i \in [n], u \in K$. By the definition of $W$, every function in $\im B$ vanishes on $W$. Since $\dim \im B=r$, we have $\dim W \le \dim V -r$. 

We now verify that $K$ and $W$ satisfy the conditions in Lemma \ref{lem_kw}. $T_i(K) \subset W$ is trivially true by the definition of $W$. So the only remaining condition to verify is $g(w, w')=0$ whenever $w, w' \in W$. Again by the definition of $W$, it suffices to show $g(T_i(u), T_j(v))=0$ for $u, v \in K$. 

Here we will use that $B$ has the maximum rank among all linear combinations of $B_i$'s. First obviously $\dim \im B=r$, so we can choose $u_1, \cdots, u_r \in U$, such that $B(u_1), \cdots, B(u_r)$ form a basis of $\im B$. We will show that $B_i(K) \subset \im B$, in other words $B_i$ sends $\ker B$ into $\im B$. Suppose not. Then there exists $u \in \ker B$ such that $B_i(u) \not \in \im B$. Then $B(u_1), \cdots, B(u_r), B_i(u)$ are $r+1$ linearly independent elements of $V^*$. Thus for sufficiently small nonzero $\varepsilon$, a small perturbation $(B+\varepsilon B_i)(u_1), \cdots, (B+\varepsilon B_i)(u_r), B_i(u)$ should still give $r+1$ linearly independent elements of $V^*$. However since $u \in \ker B$, $(B+\varepsilon B_i)(u)=\varepsilon B_i(u)$. Therefore $B+\varepsilon B_i$ has rank at least $r+1$, but it is still a linear combination of $B_i$'s, which contradicts the maximality of the rank of $B$. Hence, we always have $B_i(K) \subset \im B$. Therefore, for $u, v \in K$, we have
$g(T_i(u), T_j(v))=B_i(u)(T_j(v))$. Note that $B_i(u) \in B_i(K) \subset \im B$. Previously, we showed that every function in $\im B$ vanishes on $W$. Since $T_j(v) \in W$, we have $g(T_i(u), T_j(v))=B_i(u)(T_j(v))=0$. This immediately gives $g(w, w')=0$ for any $w, w' \in W$.

Finally, we apply Lemma \ref{lem_kw}. We have $\dim K \le \dim W$. As a consequence,
$$\dim U=\dim K +r \le \dim W+r \le \dim V.$$
Recall that $\dim U=\dim \ker O^T$ and $\dim V=\dim \ker O$, and $O$ is an $|\mathcal{F}| \times 
\left(\binom{n}{3}+n \right)$ matrix. Therefore
$$|\mathcal{F}|-\dim U = \binom{n}{3}+n-\dim V.$$
From $\dim U \le \dim V$, we immediately know that $|\mathcal{F}| \le \binom{n}{3}+n$.
\end{proof}

\section{Concluding remarks and some further discussions}
\label{sec_concluding}

The kernel and dimension arguments in this paper also lead to an alternative
elementary proof of the binary two-distance bound
$|\mathcal F|\le1+\binom n2$ for $n\ge6$, first established in
\cite{BargGlazyrinKaoLaiTsengYu}.
Our approach does not require upper bounds for equiangular line sets.
In the earlier proof, a spherical lifting transforms two inner
products into opposite values, allowing the use of such bounds.
That lifting does not directly extend to three distances: a
nonconstant affine transformation preserves three distinct inner
products, whereas unit representatives of distinct equiangular
lines have at most two.

Our main theorem determines the maximum cardinality of a $3$-distance set in the Hamming space, but leaves the
classification of extremal families unresolved.
At $n=24$, the standard construction is not the only extremal family.
Appending a fixed coordinate to the dual binary Golay code of length
$23$ gives a family of size $2048=24+\binom{24}{3}$
with distance set $\{8,12,16\}$; see
\cite[Example~9.3]{DelsarteGoethalsSeidel}.
It is therefore not isometric to the family of all vectors with
exactly one or three $1$-coordinates, whose distance set is
$\{2,4,6\}$.

For the prescribed distances $\{2,4,6\}$, the equality cases already
follow from Frankl's theorem on families with fixed
diameter~\cite{FranklFixedDiameter}.
Indeed, after a translation into the even parity class, deleting one
coordinate gives a family of diameter at most $6$ in dimension $n-1$.
Since $n+\binom n3=\sum_{j=0}^{3}\binom{n-1}{j}$,
Frankl's theorem implies that, for $n\ge9$, equality forces the
punctured family to be a Hamming ball of radius $3$.
Restoring the parity coordinate gives an isometric copy of the
standard construction.
Thus, for large $n$, the remaining issue is whether equality forces
the distance set to be $\{2,4,6\}$.

\begin{question}[Eventual uniqueness]
\label{question_eventual_uniqueness}
Does there exist an integer $n_0>24$ such that, for every $n\ge n_0$,
every binary three-distance family of size $n+\binom n3$ can be obtained
from the family of all vectors with exactly one or three
$1$-coordinates by coordinate permutations and coordinate sign
changes? If so, what is the least possible $n_0$?
\end{question}

A stronger question concerns families whose distance set differs
from $\{2,4,6\}$.
Dong, Gao, Liu, Ouyang, and Zhou~\cite{DongGaoLiuOuyangZhou}
proved that, for every fixed three-element distance set
$D\ne\{2,4,6\}$,
\[
A(n,D)\le\left(\frac34+o_D(1)\right)\binom n3.
\]
The threshold implicit in this estimate may depend on $D$.
Our theorem gives a bound uniform over $D$, but a strict inequality
below $n+\binom n3$ need not give a fixed proportional gap: several
of our estimates save only a term of order $n^2$.
We therefore ask whether such a gap can be obtained uniformly when
the distances vary with $n$.

\begin{question}[A uniform gap]
\label{question_uniform_gap}
Do there exist $\varepsilon\in(0,1)$ and an integer $n_0$ such
that, for every $n\ge n_0$ and every three-element set
$D\subseteq\{1,\dots,n\}$ with $D\ne\{2,4,6\}$,
\[
A(n,D)\le(1-\varepsilon)\left(n+\binom n3\right)?
\]
\end{question}

A positive answer to Question~\ref{question_uniform_gap}, together
with the fixed-distance equality characterization above, would
also answer Question~\ref{question_eventual_uniqueness} affirmatively.

For general $s\ge4$, our arguments reduce the range that remains in
Conjecture~\ref{conj_main}.
Suppose that $n\ge4s$. If fewer than $s$ distances occur, the
polynomial bound gives
$|\mathcal F|\le\sum_{j=0}^{s-1}\binom nj<\binom ns$.
Together with the remarks in Sections~\ref{sec_bounds} and
\ref{sec_outside}, this shows that any counterexample must have
exactly $s$ distinct distances, all even, satisfying $d_1+\cdots+d_s\le\frac{s(n+2)}2$. The coordinate pairing argument gives a further reduction.
For a family $\mathcal F\subseteq\{-1,1\}^n$ with distances in
$\{d_1,\dots,d_s\}$, where $0<d_1<\cdots<d_s\le n$, let
$P(t)=\prod_{j=1}^{s}\bigl(t-(n-2d_j)\bigr)$ and write
\[
P(\langle x,y\rangle)
=\sum_{k=0}^{s}c_k\sum_{|T|=k}\chi_T(x)\chi_T(y).
\]
Here $c_s=s!$, and we set $c_k=0$ for $k>s$.
Fix $i\in[n]$ and multiply the kernel by $1-\tau x_i y_i$.
For each $T\subseteq[n]\setminus\{i\}$ with $|T|\le s$,
the coefficients indexed by $T$ and $T\cup\{i\}$ become
\[
c_{|T|}-\tau c_{|T|+1},
\qquad
c_{|T|+1}-\tau c_{|T|}.
\]
If no two consecutive coefficients are equal and strictly positive,
then, for $\tau\in(0,1)$ sufficiently close to $1$, each pair
contains at most one strictly positive coefficient.
On $\mathcal F$, the modified kernel matrix is
$(1-\tau)2^s d_1\cdots d_s I$.
Lemma~\ref{lem_inertia} therefore gives
\[
|\mathcal F|
\le\sum_{k=0}^{s}\binom{n-1}{k}
=\sum_{\substack{0\le j\le s\\j\equiv s\pmod*2}}\binom nj.
\]
Thus this argument leaves only the cases with
$c_k=c_{k+1}>0$ for some $0\le k<s$.
Since $c_{s-1}=-(s-1)!\sum_{j=1}^{s}(n-2d_j)$,
the distance sum $d_1+\cdots+d_s=s(n+1)/2$ gives
$c_{s-1}=c_s=s!$.
This suggests studying the following extension of the central case
treated in Section~\ref{sec_difficult}.

\begin{question}[The central distance sum]
\label{question_central_sum}
For every integer $s\ge4$, does there exist an integer $N_s\ge s$,
depending only on $s$, such that the following holds?
If $n\ge N_s$ and the distinct even distances
$0<d_1<\cdots<d_s\le n$ satisfy
$d_1+\cdots+d_s=s(n+1)/2$,
must every family $\mathcal F\subseteq\{-1,1\}^n$ with distances
in $\{d_1,\dots,d_s\}$ satisfy
\[
|\mathcal F|
\le\sum_{\substack{0\le j\le s\\j\equiv s\pmod*2}}\binom nj?
\]
\end{question}

One possible approach is to extend the compatible maps and dimension
comparison in Section~\ref{sec_difficult}.
The objective is a higher-degree analogue of Lemma~\ref{lem_kw}
that gives the required dimension bound when coefficient pairing
fails.

Finally, one can ask how the threshold in
Conjecture~\ref{conj_main} depends on $s$.
Let $M_s$ be the least positive integer such that the conjectured
bound holds for every $n\ge M_s$, and set $M_s=\infty$ if no such
integer exists.
The parity construction gives $M_s\ge2s+2$, while the two-distance
theorem and the present result give $M_2=6$ and $M_3=24$.

\begin{question}[Polynomial thresholds]
\label{question_polynomial_threshold}
Does there exist an absolute constant $C>0$ such that
\[
M_s\le s^C\qquad\text{for every integer }s\ge2?
\]
\end{question}

\appendix
\section{Proof of Lemma \ref{lem_la}}\label{app_1}
\begin{proof}
We first prove that for $P, Q \in \mathcal{Q}$, $Qx=0$ implies $QPx=0$. Fix $Q \in \mathcal{Q}$. We can choose the coordinate system such that 
$$Q=\textup{diag}(q_1, \cdots, q_r, 0, \cdots, 0)$$
for $q_i \neq 0$. To justify this change of coordinates, let $S=(s_{ab})$ be an
orthogonal matrix. For each $Q\in\mathcal Q$, write $Q'=S^TQS$, and define $R'_{Q'}$ by
\[
(R'_{Q'})_{i,(j,k,\ell)}
=
\sum_{a,b,c,d=1}^{n}
s_{ai}s_{bj}s_{ck}s_{d\ell}
(R_Q)_{a,(b,c,d)}.
\]
Since $S$ is invertible, $Q'$ uniquely determines $Q$. Thus this
defines a linear map on the transformed space
$\mathcal Q'=\{S^TQS:Q\in\mathcal Q\}$. The entries of $R'_{Q'}$ are invariant under permutations of
$i,j,k,\ell$, because the entries of $R_Q$ are invariant under
permutations of $a,b,c,d$. Moreover, for $P'=S^TPS$, the
orthogonality of $S$ gives
\[
R'_{Q'}(R'_{P'})^T
=
S^T(R_QR_P^T)S
=
S^TQPS
=
Q'P'.
\]
Hence both assumptions of the lemma are preserved under the change
of coordinates. We henceforth work in the new coordinates and omit
the primes. From the assumption, $R_QR_Q^T=Q^2$. Note that the last $n-r$ diagonal entries of $Q^2$ are $0$, and the $i$-th diagonal entry of $R_Q R_Q^T$ is the sum of squares of entries in row $i$ of $R_Q$. Therefore the last $n-r$ rows of $R_Q$ are zero vectors. By the assumption that $(R_Q)_{i,(j,k,\ell)}$ is invariant under permutations of $i,j,k, \ell$, we know an entry of $R_Q$ can be non-zero only when all $i,j,k,\ell$ are at most $r$. Also, the first $r$ rows of $R_Q$, when  divided by $q_1, \cdots, q_r$ respectively, are mutually orthogonal unit vectors. We call these vectors $e_1, \cdots, e_r$. 

Now fix $P\in\mathcal Q$, and let $v$ be the $i$-th row of $R_P$.
Using $R_PR_Q^T=PQ$, we have, for $j=1,\ldots,r$,
\[
\langle v,e_j\rangle
=
\frac{(R_PR_Q^T)_{ij}}{q_j}
=
\frac{(PQ)_{ij}}{q_j}
=
P_{ij}.
\]

Let $v_{\le r}$ be the restriction of $v$ to the coordinates indexed
by triples $(j,k,\ell)$ with $j,k,\ell\le r$, with all other
coordinates set to zero. By the permutation symmetry, each $e_j$
is supported on these coordinates. Hence
\[
\langle v,e_j\rangle=\langle v_{\le r},e_j\rangle
\qquad\text{for }j=1,\ldots,r.
\]
Since $e_1,\ldots,e_r$ are orthonormal, Bessel's inequality gives
\begin{equation}\label{eq_RP}
\sum_{j=1}^{r}P_{ij}^{2}
=
\sum_{j=1}^{r}\langle v_{\le r},e_j\rangle^{2}
\le
\|v_{\le r}\|_2^{2}
=
\sum_{j,k,\ell\le r}
(R_P)_{i,(j,k,\ell)}^{2}.
\end{equation}
We let $$b=\sum_{i=1}^r \sum_{j=r+1}^n P_{ij}^2 \ge 0.$$
Our goal is to show that $b=0$ by estimating the sum of the squares of the entries of $R_P$ in two different ways. First, summing \eqref{eq_RP} for $i>r$ gives 
$$\sum_{i>r, j,k, \ell \le r} (R_P)_{i,(j,k,\ell)}^2 \ge \sum_{i>r}\sum_{j \le r} P_{ij}^2=b.$$
Then the permutation symmetry of $R_P$ immediately gives
$$\sum_{i \le r, \textup{exactly~one~of~}j,k,\ell>r} (R_P)^2_{i,(j,k,\ell)} \ge 3b.$$
On the other hand, since $R_P R_P^T=P^2$, the squared length of the $i$-th row of $R_P$ is equal to $\sum_{j=1}^n P_{ij}^2$. \eqref{eq_RP} already says that the contribution for $(j,k,\ell)$ with all three indices at most $r$ is already at least $\sum_{j=1}^r P_{ij}^2$. Subtracting these, we get
$$\sum_{i \le r, \textup{exactly~one~of~}j,k,\ell>r} (R_P)^2_{i,(j,k,\ell)} \le \sum_{i \le r} \left(\sum_{j=1}^n P_{ij}^2-\sum_{j \le r} P_{ij}^2\right)=\sum_{i\le r, j>r} P_{ij}^2 = b.$$
Therefore $3b \le b$ and since $b \ge 0$, we have $b=0$. Consequently, $P_{ij}=0$ whenever $i \le r$ and $j>r$. Now if $Qx=0$, we know the first $r$ coordinates of $x$ must vanish, but then this would imply that the first $r$ coordinates of $Px$ vanish as well by this property of $P$. Then $QPx=0$ is automatically true.

We choose $x \in \mathbb{R}^n$ such that the subspace $\{P x: P \in \mathcal{Q}\}$ has the largest possible dimension, say $d$. We choose $P_1 x, \cdots, P_d x$ as a basis of this subspace. We claim that $Qx=0$ for $Q \in \mathcal{Q}$ forces $Q=0$. Suppose instead we have $Q \neq 0$ and $Qx=0$. By the previous observation, we also have $QP_ix=0$ for $i=1, \cdots, d$. Choose $y \in \mathbb{R}^n$ so that $Qy \neq 0$. Then, since $Q$ is symmetric,
$$\langle Qy, P_i x\rangle=\langle y, QP_i x\rangle =0.$$
So $Qy$ is a non-zero vector orthogonal to all $P_i x$. Hence $P_1 x, \cdots, P_d x, Qy$ are linearly independent. We use the ``small perturbation'' trick used in the proof of Lemma \ref{lem_uv} once again, and have that for sufficiently small nonzero $t$,
$$P_1(x+ty), \cdots, P_d(x+ty), Qy$$
are still linearly independent. Since $Qx=0$, we have $Q(x+ty)=tQy$. So we know $$P_1(x+ty), \cdots, P_d(x+ty), Q(x+ty)$$
are linearly independent. This gives $x+ty$ such that $\{P(x+ty): P \in \mathcal{Q}\}$ has dimension at least $d+1$, contradicting the maximality assumption.

Therefore $Qx=0$ implies $Q=0$, meaning that $Q \rightarrow Qx$ is injective. The image of this mapping is contained in $\sum_{Q \in \mathcal{Q}} \im Q$, so 
$$\dim \mathcal{Q} \le \dim \left(\sum_{Q \in \mathcal{Q}} \im Q\right),$$
which completes the proof of Lemma \ref{lem_la}.
\end{proof}

\end{document}